\documentclass[reqno]{amsart}
\usepackage{amssymb}
\usepackage{ifpdf}
\ifpdf
 \usepackage[hyperindex,pagebackref]{hyperref}
\else
 \expandafter\ifx\csname dvipdfm\endcsname\relax
 \usepackage[hypertex,hyperindex,pagebackref]{hyperref}
 \else
 \usepackage[dvipdfm,hyperindex,pagebackref]{hyperref}
 \fi
\fi
\allowdisplaybreaks[4]
\theoremstyle{plain}
\newtheorem{theorem}{Theorem}[section]
\theoremstyle{remark}
\newtheorem{remark}{Remark}
\DeclareMathOperator{\td}{d}
\numberwithin{equation}{section}
\date{Completed on 9 February 2011}
\date{}

\begin{document}

\title[A completely monotonic function involving ${p}$-psi function]
{Complete monotonicity of a function involving the $\boldsymbol{p}$-psi function and alternative proofs}

\author[V. Krasniqi]{Valmir Krasniqi}
\address[Krasniqi]{Department of Mathematics, University of Prishtina, Prishtin\"{e} 10000, Republic of Kosova}
\email{\href{mailto: V. Krasniqi <vali.99@hotmail.com>}{vali.99@hotmail.com}}

\author[F. Qi]{Feng Qi}
\address[Qi]{School of Mathematics and Informatics, Henan Polytechnic University, Jiaozuo City, Henan Province, 454010, China; Department of Mathematics, School of Science, Tianjin Polytechnic University, Tianjin City, 300387, China}
\email{\href{mailto: F. Qi <qifeng618@gmail.com>}{qifeng618@gmail.com}, \href{mailto: F. Qi <qifeng618@hotmail.com>}{qifeng618@hotmail.com}, \href{mailto: F. Qi <qifeng618@qq.com>}{qifeng618@qq.com}}
\urladdr{\url{http://qifeng618.wordpress.com}}



\begin{abstract}
In the paper the authors alternatively prove that the function $x^\alpha\big[\ln\frac{px}{x+p+1}-\psi_p(x)\big]$ is completely monotonic on $(0,\infty)$ if and only if $\alpha \le 1$, where $p\in\mathbb{N}$ and $\psi_p(x)$ is the $p$-analogue of the classical psi function $\psi(x)$. This generalizes a known result.
\end{abstract}

\maketitle

\section{Introduction}

Recall from~\cite[Chapter~XIII]{mpf-1993}, \cite[Chapter~1]{Schilling-Song-Vondracek-2010} and~\cite[Chapter~IV]{widder} that a function $f$ is said to be completely monotonic on an interval $I$ if $f$ has derivatives of all orders on $I$
and satisfies
\begin{equation}\label{CM-dfn}
0\le(-1)^{n}f^{(n)}(x)<\infty
\end{equation}
for $x\in I$ and $n\ge0$.
The celebrated Bernstein-Widder's Theorem (see~\cite[p.~3, Theorem~1.4]{Schilling-Song-Vondracek-2010} or~\cite[p.~161, Theorem~12b]{widder}) characterizes that a necessary and sufficient condition that $f(x)$ should be completely monotonic for $0<x<\infty$ is that
\begin{equation} \label{berstein-1}
f(x)=\int_0^\infty e^{-xt}\td\alpha(t),
\end{equation}
where $\alpha(t)$ is non-decreasing and the integral converges for $0<x<\infty$. This expresses that a completely monotonic function $f$ on $[0,\infty)$ is a Laplace transform of the measure $\alpha$.
\par
It is common knowledge that the classical Euler's gamma function $\Gamma(x)$ may be defined for $x>0$ by
\begin{equation*}
\Gamma(x)=\int_0^\infty t^{x-1}e^{-t}\td t.
\end{equation*}
The logarithmic derivative of $\Gamma(x)$, denoted by $\psi(x)=\frac{\Gamma'(x)}{\Gamma(x)}$, is called psi function or digamma function.
\par
An alternative definition of the gamma function $\Gamma(x)$ is
\begin{equation}\label{eqgammap2}
\Gamma(x) =\lim_{p\rightarrow\infty}\Gamma_p(x),
\end{equation}
where
\begin{equation}\label{eqgammap1}
\Gamma_p(x)=\frac{p!p^x}{x(x+1)\dotsm(x+p)}=\frac{p^x}{x(1+{x}/{1})\dotsm(1+{x}/{p})}
\end{equation}
for $x>0$ and $p\in\mathbb{N}$, the set of all positive integers. See~\cite[p.~250]{Apostol-Number}. The $p$-analogue of the psi function $\psi(x)$ is defined as the
logarithmic derivative of the $\Gamma_p$ function, that is,
\begin{equation}\label{psi_p1}
\psi_p(x)=\frac{\td}{\td x} \ln \Gamma_p(x)=\frac{\Gamma'_p(x)}{\Gamma_p(x)}.
\end{equation}
The function $\psi_p$ has the following properties:
\begin{enumerate}
\item
It has the following representations
\begin{equation}\label{psi_series1}
\psi_p(x)=\ln p-\sum_{k=0}^{p}\frac{1}{x+k} =\ln p -\int_{0}^{\infty}\frac{1-e^{-(p+1)t}}{1-e^{-t}}e^{-xt}\td t.
\end{equation}
\item
It is increasing on $(0,\infty)$ and $\psi'_p$ is completely monotonic on $(0,\infty)$.
\end{enumerate}
The very right hand side of the formula~\eqref{psi_series1} corrects errors appeared in~\cite[p.~374, Lemma~5]{MC-Krashiqi-Mansour-Shabani} and~\cite[p.~29, Lemma~2.3]{VA}.
\par
In~\cite[pp.~374\nobreakdash--375, Theorem~1]{psi-alzer}, it was proved that the function
\begin{equation}\label{theta-alpha}
\theta_\alpha(x)=x^\alpha[\ln x-\psi(x)]
\end{equation}
is completely monotonic on $(0,\infty)$ if and only if $\alpha\le1$. For the history, background, applications and alternative proofs of this conclusion, please refer to~\cite{theta-new-proof.tex-BKMS}, \cite[p.~8, Section~1.6.6]{bounds-two-gammas.tex} and closely related references therein.
\par
The aim of this paper is to generalize~\cite[pp.~374\nobreakdash--375, Theorem~1]{psi-alzer} and~\cite[p.~105, Theorem~1]{theta-new-proof.tex-BKMS} to the case of the $p$-analogue of the psi function $\psi(x)$ as follows.

\begin{theorem}\label{vali-Qi-thm}
The function
\begin{equation}
\theta_{p,\alpha}(x)=x^\alpha \biggl[\ln\frac{px}{x+p+1}-\psi_p(x)\biggr]
\end{equation}
for $p\in\mathbb{N}$ is completely monotonic on $(0,\infty)$ if and only if $\alpha \le 1$.
\end{theorem}

\begin{remark}
Letting $p\to\infty$ in~Theorem~\ref{vali-Qi-thm}, we obtain~\cite[pp.~374\nobreakdash--375, Theorem~1]{psi-alzer} and~\cite[p.~105, Theorem~1]{theta-new-proof.tex-BKMS}.
\end{remark}

\section{Proofs of Theorem~\ref{vali-Qi-thm}}

\begin{proof}[First Proof]
From the identity~\eqref{psi_series1} and the integral expression
\begin{equation}\label{log}
\ln\frac{b}a=\int_{0}^{\infty}\frac{e^{-at}-e^{-bt}}{t}\td t
\end{equation}
in~\cite[p.~230, 5.1.32]{abram}, we obtain
\begin{equation}\label{vali-qi-theta-int}
\theta_{p,1}(x)=x\int_{0}^{\infty}\bigl[1-e^{-(p+1)t}\bigr]\varphi (t)e^{-xt}\td t,
\end{equation}
where
\begin{equation}
\varphi (t)=\frac{1}{1-e^{-t}}-\frac{1}{t}.
\end{equation}
The function $\varphi(t)$ is increasing on $(0,\infty)$ with
\begin{equation}\label{varphi-2-limits}
\lim_{t\to 0}\varphi (t)=\frac{1}{2}\quad \text{and} \quad \lim_{t\to \infty}\varphi (t)=1.
\end{equation} 
See~\cite{best-constant-one-simple.tex, emv-log-convex-simple.tex, remiander-Guo-Qi-Tamsui, best-constant-one.tex, Cheung-Qi-Rev.tex, Extended-Binet-remiander-comp.tex, best-constant-one-simple-real.tex} and related references therein.
Therefore, for $x>0$ and $n\in\mathbb{N}$, we have
\begin{align*}
(-1)^n\theta_{p,1} ^{(n)}(x) & =x(-1)^n\frac{\td^n}{\td x^n}
\int_{0}^{\infty}\bigl[1-e^{-(p+1)t}\bigr]\varphi(t)e^{-xt}\td t \\
 &\quad-(-1)^{n-1}n\frac{\td^{n-1}}{\td x^{n-1}}\int_{0}^{\infty}\bigl[1-e^{-(p+1)t}\bigr]\varphi (t)e^{-xt}\td t \\
 &=x\int_{0}^{\infty}t^{n}\varphi (t)\bigl[1-e^{-(p+1)t}\bigr]e^{-xt}\td t\\
 &\quad-n\int_{0}^{\infty}t^{n-1}\varphi (t)\bigl[1-e^{-(p+1)t}\bigr]e^{-xt}\td t \\
 &=\int_{0}^{n/x}t^{n-1}\bigl[1-e^{-(p+1)t}\bigr]\varphi (t)(tx-n)e^{-xt}\td t \\
 &\quad+\int_{n/x}^{\infty}t^{n-1}\bigl[1-e^{-(p+1)t}\bigr]\varphi (t)(tx-n)e^{-xt}\td t\\
&>\varphi \biggl(\frac{n}{x}\biggl)\int_{0}^{n/x}t^{n-1}\bigl[1-e^{-(p+1)t}\bigr](tx-n)e^{-xt}\td t \\
&\quad+\varphi \biggl(\frac{n}{x}\biggl)\int_{n/x}^{\infty}t^{n-1}\bigl[1-e^{-(p+1)t}\bigr](tx-n)e^{-xt}\td t \\
&=\varphi \biggl(\frac{n}{x}\biggl)\int_{0}^{\infty}t^{n-1}\bigl[1-e^{-(p+1)t}\bigr](tx-n)e^{-xt}\td t \\
&=\varphi \biggl(\frac{n}{x}\biggl)\biggl[x\int_{0}^{\infty}t^{n}\bigl[1-e^{-(p+1)t}\bigr]e^{-xt}\td t\\
&\quad-n\int_{0}^{\infty}t^{n-1}\bigl[1-e^{-(p+1)t}\bigr]e^{-xt}\td t\biggl]\\
&=\varphi \biggl(\frac{n}{x}\biggl)\biggl[x\int_{0}^{\infty}t^{n}e^{-xt}\td t
-x\int_{0}^{\infty}t^{n}e^{-(x+p+1)t}\td t\\
&\quad-n\int_{0}^{\infty}t^{n-1}e^{-xt}\td t+n\int_{0}^{\infty}t^{n-1}e^{-(x+p+1)t}\td t\biggl] \\
&=\varphi \biggl(\frac{n}{x}\biggl)\biggr[ x\frac{n!}{x^{n+1}}- x\frac{n!}{(x+p+1)^{n+1}}- n\frac{(n-1)!}{x^{n}}+ n\frac{(n-1)!}{(x+p+1)^{n}}\biggr]\\
&=\varphi \biggl(\frac{n}{x}\biggl)n!\biggr[ \frac1{x^{n}}- \frac{x}{(x+p+1)^{n+1}}- \frac{1}{x^{n}}+ \frac{1}{(x+p+1)^{n}}\biggr]\\
&=\varphi \biggl(\frac{n}{x}\biggl)\frac{n!}{(x+p+1)^{n}}\biggr(1- \frac{x}{x+p+1}\biggr)\\
&=\varphi \biggl(\frac{n}{x}\biggl)\frac{n!(p+1)}{(x+p+1)^{n+1}}\\
&>0,
\end{align*}
where we used the formula
\begin{equation}\label{fracint}
\frac1{x^\omega}=\frac1{\Gamma(\omega)}\int_0^\infty t^{\omega-1}e^{-xt}\td t
\end{equation}
for real numbers $x>0$ and $\omega>0$, see~\cite[p.~255, 6.1.1]{abram}. So we obtain that the function $\theta_{p,1}(x)$ is completely monotonic on $(0,\infty)$.
\par
Since
\begin{equation*}
(-1)^n[u(x) v(x)]^{(n)}=\sum_{i=0}^{n}\binom{n}{i} \bigl[(-1)^{i}u^{(i)}(x)\bigr]\bigl[(-1)^{n-i}v^{(n-i)}(x)\bigr],
\end{equation*}
the product of any two completely monotonic function is also completely monotonic on their common domain. On the other hand, the function $x^{\alpha-1}$ for $\alpha <1$ is clearly completely monotonic on $(0,\infty)$. Consequently the function
\begin{equation*}
\theta_{p,\alpha}(x)=x^{\alpha-1}\theta_{p,1}(x)
\end{equation*}
for $\alpha\le1$ is completely monotonic on $(0,\infty)$.
\par
Conversely,  if $\theta_{p,\alpha}(x)$ is completely monotonic on $(0,\infty)$, then
\begin{equation*}
\frac{\td\theta_{p,\alpha}(x)}{\td x}=x^{\alpha-1}\biggr\{\alpha\biggl[\ln\frac{px}{x+p+1}
-\psi_p(x)\biggl] +\frac{p}{x+p+1}-x\psi'_p(x)\biggr\}\le0
\end{equation*}
for $x>0$, equivalently,
\begin{equation*}
\alpha\le\frac{x\psi'_p(x)-\frac{p}{x+p+1}}{\ln\frac{px}{x+p+1}-\psi_p(x)}.
\end{equation*}
Employing L'H\^ospital's rule and~\eqref{psi_series1} results in
\begin{align*}
\lim_{x\to\infty}\frac{x\psi'_p(x)-\frac{p}{x+p+1}}{\ln\frac{px}{x+p+1}-\psi_p(x)}
&=\lim_{x\to\infty}\frac{x\psi_p''(x)+\psi_p'(x)+\frac{p}{(x+p+1)^2}}{\frac1x-\frac1{x+p+1}-\psi_p'(x)}\\
&=\lim_{x\to\infty}\frac{\frac{p}{(x+p+1)^2}-x\sum_{k=0}^{p}\frac{2}{(x+k)^3}+\sum_{k=0}^{p}\frac{1}{(x+k)^2}} {\frac1x-\frac1{x+p+1}-\sum_{k=0}^{p}\frac{1}{(x+k)^2}}\\
&=1,
\end{align*}
so it is necessary that $\alpha\le1$. The proof is complete.
\end{proof}

\begin{proof}[Second Proof]
From~\eqref{vali-qi-theta-int} and by integration by part lead to
\begin{align*}
\theta_{p,1}(x)&=-\int_{0}^{\infty}\bigl[1-e^{-(p+1)t}\bigr]\varphi (t)\frac{\td e^{-xt}}{\td t}\td t\\
&=\int_{0}^{\infty}\bigl\{\bigl[1-e^{-(p+1)t}\bigr]\varphi (t)\bigr\}'e^{-xt}\td t
-\bigl\{\bigl[1-e^{-(p+1)t}\bigr]\varphi (t)e^{-xt}\bigr\}\big\vert_{t=0}^{t=\infty}\\
&=\int_{0}^{\infty}\bigl\{\bigl[1-e^{-(p+1)t}\bigr]\varphi'(t)+(p+1)e^{-(p+1)t}\varphi(t)\bigr\}e^{-xt}\td t.
\end{align*}
Therefore, for showing that the function $\theta_{p,1}(x)$ is completely monotonic on $(0,\infty)$ for all $p\in\mathbb{N}$, it suffices to prove that the function 
\begin{equation}\label{varphi-exp-positivity}
\bigl[1-e^{-(p+1)t}\bigr]\varphi'(t)+(p+1)e^{-(p+1)t}\varphi(t)
\end{equation}
is positive. Since the function $\varphi(t)$ is increasing on $(0,\infty)$, the derivative $\varphi'(t)$ is positive on $(0,\infty)$. Further considering the limits in~\eqref{varphi-2-limits}, the positivity of $\varphi(t)$ follows. As a result, the function~\eqref{varphi-exp-positivity} is positive.
\par
The rest of the proof is the same as the first proof.
\end{proof}

\subsection*{Acknowledgements}
The authors are indebted to anonymous referees for pointing out some errors appeared in the first version of this manuscript.

\end{document}